\documentclass[12pt]{article}
\usepackage{amsmath}
\usepackage{amssymb}
\usepackage{vatola}

\def\q{\quad}

\def\qtq#1{\q\t{#1}\q}

\def\({\left(}
\def\){\right)}
\def\[{\left[}
\def\]{\right]}
\def\mod#1{\ (\text{\rm mod}\ #1)}
\def\t{\text}
\def\f{\frac}

\def\e{\equiv}

\def\b{\binom}

\def\cos{\t{\rm cos}}
\def\sin{\t{\rm sin}}

\def\sls#1#2{(\f{#1}{#2})}

\def\Ls#1#2{\Big(\f{#1}{#2}\Big)}
\let \pro=\proclaim
\let \endpro=\endproclaim

\begin{document}
 \centerline {\bf
On the properties of invariant functions}
\par\q\newline
\centerline{Zhi-Hong Sun}\newline \centerline{School of Mathematics
and Statistics} \centerline{Huaiyin Normal University}
\centerline{Huaian, Jiangsu 223300, P.R. China} \centerline{Email:
zhsun@hytc.edu.cn} \centerline{Homepage:
http://maths.hytc.edu.cn/szh1.htm}
 \abstract{If $f(x,y)$ is a real function satisfying $y>0$ and
 $\sum_{r=0}^{n-1}f(x+ry,ny)=f(x,y)$ for $n=1,2,3,\ldots$, we say
 that $f(x,y)$ is an invariant function. Many special functions
 including Bernoulli polynomials, Gamma function and Hurwitz zeta
 function are related to invariant functions.
 In this paper we systematically investigate the properties of invariant
 functions.
 \par\q
 }
 \endabstract

\section*{1. Introduction}
\par Let $\Bbb Z$ and $\Bbb Z^+$ be the set of integers and the set of positive
integers, respectively. For $a\in\Bbb Z$ and $m\in\Bbb Z^+$ let
$a(m)=a+m\Bbb Z=\{a+km\bigm| k\in\Bbb Z\}.$ If $$a_1(n_1)\cup
a_2(n_2)\cup\cdots\cup a_k(n_k)=\Bbb Z\qtq{and}a_i(n_i)\cap
a_j(n_j)=\emptyset\qtq{for any}i\not=j,$$ we say that
$\{a_1(n_1),\ldots,a_k(n_k)\}$ is a disjoint covering system. In
[6,7,8], Z.W. Sun showed that if $\{a_1(n_1),\ldots,a_k(n_k)\}$ is a
disjoint covering system and $F(x,y)$ satisfies
$$\sum_{r=0}^{n-1}F\Big(\f{x+r}n,ny\Big)=F(x,y)\qtq{for}n=1,2,3,\ldots,\tag 1.1$$
then $\sum_{s=1}^kF(\f{x+a_s}{n_s},{n_s}y)=F(x,y).$ In [7], Z.W. Sun
also gave some examples of $F(x,y)$ satisfies (1.1). \par  If
$f(x,y)$ is a real function with $y>0$ and
$$\sum_{r=0}^{n-1}f(x+ry,ny)=f(x,y)\qtq{for}n=1,2,3,\ldots,\tag 1.2$$
 we say that $f(x,y)$
is an invariant function. This is because
$$f(x+y,y)-f(x,y)=\sum_{r=0}^{n-1}f\big(x+\f yn+\f
rny,y\big)-\sum_{r=0}^{n-1}f\big(x+\f rny,y\big)=f\big(x+\f yn,\f
yn\big)-f\big(x,\f yn\big)$$ for any $n\in\Bbb Z^+$ and so
$$f(x+y,y)-f(x,y)=\lim_{n\rightarrow +\infty}\Big(f\big(x+\f yn,\f
yn\big)-f\big(x,\f yn\big)\Big)=\lim_{a\rightarrow
0^+}(f(x+a,a)-f(x,a)).\tag 1.3$$ If the invariant function $f(x,y)$
is an integrable function of $x$ at any closed intervals, we say
that $f(x,y)$ is an integrable invariant function. The set of
integrable invariant functions is denoted by $I$. If $F(x,y)$ is a
real function satisfying (1.1) and $f(x,y)=F(\f xy,y)$, then clearly
$f$ is an invariant function. By the above result due to Z.W. Sun,
if $\{a_1(n_1),\ldots,a_k(n_k)\}$ is a disjoint covering system and
$f$ is an invariant function, then
$$\sum_{s=1}^kf(x+a_sy,n_sy)=f(x,y)\qtq{and so}\sum_{s=1}^kf(a_s,n_s)=f(0,1).
\tag 1.4$$

\par Let $[x]$ be the greatest integer not exceeding $x$. By
the Hermite identity $\sum_{r=0}^{n-1}[x+\f rn]=[nx]$, $[\f xy]$ is
an invariant function. The
 Bernoulli numbers $\{B_n\}$ and Bernoulli polynomials
$\{B_n(x)\}$ are given by
$$B_0=1,\q\sum_{k=0}^{n-1}\b nkB_k=0\q(n\ge 2)\qtq{and}
B_n(x)=\sum_{k=0}^n\b nkB_kx^{n-k}\ (n\ge 0).$$ Raabe's theorem
states (see [4]) that $$\sum_{r=0}^{n-1}B_m\Big(x+\f
rn\Big)=n^{1-m}B_m(nx)\q (m\ge 0,n\ge 1).$$ Thus, $y^{m-1}B_m(\f
xy)$ is an invariant function. More generally, if $a$ is a real
number and $F(x)$ is a real function satisfying
$$\sum_{r=0}^{n-1}F\Big(x+\f
rn\Big)=n^{-a}F(nx)\qtq{for any} n\in\Bbb Z^+,\tag 1.5$$ one can
easily check that $F(\f xy)y^a$ is an invariant function. We note
that the functional equation (1.5) have been investigated by several
mathematicians including Bass, Kubert, Milnor and Walum. See
[2],[3],[5] and [9].
\par In this paper, we systematically investigate the properties of
invariant functions. In Section 2, we point out basic properties and
more examples of invariant functions. In Section 3, we prove some
interesting results for integrable invariant functions. In
particular, if $f(x,y)$ is an integrable invariant function and
$\f{\partial f} {\partial y}$ exists, then
$$f(x,y)=\int_x^{x-y}\f{\partial}{\partial y}f(t,y)\; dt;$$
 if $g$ and $h$ are
integrable invariant functions and
$$g*h(x,y)=\int_0^xg(t,y)h(x-t,y)dt+\int_x^yg(t,y)h(x+y-t,y)dt,$$
then $g*h$ is also an integrable invariant function and
$$\int_0^yg*h(x,y)dx=\int_0^yg(x,y)dx\cdot \int_0^yh(x,y)dx.$$
For $m,n\in\Bbb Z^+$ we have
$$B_{m+n}(x)=-\binom{m+n}m\Big(\int_0^1B_m(x-t)B_n(t)dt
+m\int_x^1(x-t)^{m-1}B_n(t)dt\Big)$$ and so
$$\f{-y^{m-1}B_m\sls xy}{m!}*\f{-y^{n-1}B_n\sls xy}{n!}=
\f{-y^{m+n-1}B_{m+n}\sls xy}{(m+n)!}.$$

\section*{2. Basic properties and examples of invariant functions}
\par Let $\Bbb R$, $\Bbb R^+$ and $\Bbb C$
 be the set of real numbers, the set of positive real numbers and the set of
 complex numbers, respectively. From the definition of invariant function one can easily
prove the following three propositions.
 \pro{Proposition 2.1} Suppose that $f(x,y)$ is an
invariant function,
  $a,b,c\in\Bbb R$, $c>0$ and $F(x,y)=af(b+cx,cy)$, then
 $F(x,y)$ is also an
invariant function.
\endpro
\pro{Proposition 2.2} If $f$ is an invariant function and
$\f{\partial f}{\partial x}$ exists, then $\f{\partial f}{\partial
x}$ is also an invariant function.\endpro

\pro{Proposition 2.3} Let $F(x,y)$ be the mapping from $\Bbb R\times
\Bbb R^+$ to $\Bbb C$, and
 $F(x,y)=f(x,y)+ig(x,y)$, where $f(x,y)$ is
the real part of $F(x,y)$. If $\sum_{r=0}^{n-1}F(x+ry,ny)=F(x,y)$
for any $n\in\Bbb Z^+$, then $f(x,y)$ and $g(x,y)$ are invariant
functions.

\pro{Proposition 2.4} Suppose that $f(x,y)$ is an invariant function
and $m,n\in\Bbb Z^+$. Then
$$\sum_{r=0}^{n-1}f(x+rmy,ny)=\sum_{r=0}^{m-1}f(x+rny,my).$$
\endpro
{\it Proof.} By (1.2),
$$\align \sum_{r=0}^{n-1}f(x+rmy,ny)
&=\sum_{r=0}^{n-1}\sum_{k=0}^{m-1}f(x+rmy+k(ny),m(ny))
\\&=\sum_{k=0}^{m-1}\sum_{r=0}^{n-1}f(x+kny+r(my),n(my))
=\sum_{k=0}^{m-1}f(x+kny,my).\endalign$$
 \pro{Proposition 2.5} If $f$ is an
invariant function and $F(x,y)=f(y-x,y)$, then $F$ is also an
invariant function. \endpro {\it Proof}. For $n\in\Bbb Z^+$,
$$\align \sum_{r=0}^{n-1}
F(x+ry,ny)&=\sum_{r=0}^{n-1}f(ny-(x+ry),ny)=\sum_{r
=0}^{n-1}f(y-x+(n-1-r)y,ny)
\\&=\sum_{k=0}^{n-1}f(y-x+ky,ny)=f(y-x,y)=F(x,y).\endalign$$
This proves the proposition.
\par\q
\par
For $x\in\Bbb R$ let $\{x\}$ be the fractional part of $x$. That is,
$\{x\}=x-[x]$.
 \pro{Proposition 2.6} Suppose that $f$ is an invariant function
 and $t\in\Bbb R$. Let
$$f_1(x,y)=f\Big(y\Big\{\f{t+x}y\Big\},y\Big)
\qtq{and}f_2(x,y)=f\Big(y\Big\{\f{t-x}y\Big\},y\Big).$$
 Then both
$f_1$ and $f_2$ are invariant functions.
\endpro
{\it Proof.} For $n\in\Bbb Z^+$ we have
$$\align &\sum_{r=0}^{n-1}f_1(x+ry,ny)
\\&=\sum_{r=0}^{n-1}f\Big(ny\Big\{\f{t+x+ry}{ny}\Big\},ny\Big)
=\sum_{r=0}^{n-1}f\Big(ny\Big\{\f{r+[\f {t+x}y]+\{\f{t+
x}y\}}n\Big\},ny\Big)
\\&=\sum_{r=0}^{n-1}f\Big(ny\Big\{\f{r+\{\f{t+x}y\}}n\Big\},ny\Big)
=\sum_{r=0}^{n-1}f\Big(y\big\{\f{t+x}y\big\}+ry,ny\Big)
=f\Big(y\Big\{\f{t+x}y\Big\},y\Big).
\endalign$$
This shows that $f_1(x,y)$ is an invariant function. Since
$f_2(x,y)=f_1(y-x,y)$, $f_2(x,y)$ is also an invariant function by
Proposition 2.5. This completes the proof.

 \pro{Proposition 2.7} Let $h(x)$ be a real function, and let
 $$\f 1y\sum_{k=1}^{\infty}h\Ls ky\cos 2\pi\f {kx}y=f(x,y)
 \qtq{and}\f 1y\sum_{k=1}^{\infty}h\Ls ky\sin 2\pi\f {kx}y=g(x,y).$$
Then $f$ and $g$ are invariant functions.
\endpro
{\it Proof.} Suppose
$$\f 1y\sum_{k=1}^{\infty}h\Ls ky\t{e}^{2\pi i \f{kx}y}=F(x,y).$$
For $n\in\Bbb Z^+$ we see that
$$\align \sum_{r=0}^{n-1}F(x+ry,ny)
&=\sum_{r=0}^{n-1}\f 1{ny}\sum_{k=1}^{\infty}h\Ls k{ny}\t{e}^{2\pi i
\f{k(x+ry)}{ny}}
\\&=\f 1{ny}\sum_{k=1}^{\infty}h\Ls k{ny}\t{e}^{2\pi
i\f{kx}{ny}}\sum_{r=0}^{n-1}\t{e}^{2\pi i\f{kr}n}
\\&=\f 1{ny}\sum_{m=1}^{\infty}h\Ls my\t{e}^{2\pi i\f{mx}y}\cdot n
=F(x,y).\endalign$$ By Euler's formula
$\t{e}^{i\theta}=\cos\theta+i\sin\theta$, we have
$F(x,y)=f(x,y)+ig(x,y)$. Thus the result follows from Proposition
2.3.

 \pro{Proposition 2.8} Let $h(x)$ be a real
function and $\sum_{k=0}^{\infty}h(x+ky)=f(x,y)$ for $x\in\Bbb R$
and $y\in\Bbb R^+$. Then $f$ is an invariant function.
\endpro
{\it Proof.} For $n\in\Bbb Z^+$,
$$\sum_{r=0}^{n-1}f(x+ry,ny)=\sum_{r=0}^{n-1}
\sum_{k=0}^{\infty}h(x+ry+kny)
=\sum_{m=0}^{\infty}h(x+my)=f(x,y).$$
\par{\bf Remark 2.1} Let $h(x)$ be a real
function and $\sum_{k=-\infty}^{\infty}h(xy+ky)=F(x,y)$ for
$x\in\Bbb R$ and $y\in\Bbb R^+$. In [6,7], Z.W. Sun stated that
$F(x,y)$ satisfies (1.1).

\par\q
\par Now we list some typical examples of invariant functions.
\par{\bf Example 2.1.} $\f 1y\in I$.
\par{\bf Example 2.2.} $y^{m-1}B_m(\f xy)\in I$ $(m\in\Bbb Z^+)$.
\par{\bf Example 2.3.} $[\f xy]\in I$ and $\{\f xy\}-\f 12\in I$.
\par Since $[\f xy]\in I$ and $\f xy-\f 12=B_1\sls xy\in I$
we have $\{\f xy\}-\f 12\in I$.
 \par{\bf Example 2.4.} For $a\in\Bbb R$ let
 $$f(x,y)=\cases 1&\t{if $\f{a-x}y\in\Bbb Z,$}
\\0&\t{if $\f{a-x}y\notin\Bbb Z.$}\endcases$$
Then $f\in I$. This is because $\f{a-x}y\in\Bbb Z$ implies that
there is a unique $r\in\{0,1,\ldots,n-1\}$ such that $\f{a-x}y\e
r\mod n$ for $n\in\Bbb Z^+$.
\par{\bf Example 2.5.} For $a>0$ and $a\not=1$, we have
$\f{a^x}{a^y-1}\in I$.
\par This is because
$$\sum_{r=0}^{n-1}\f{a^{x+ry}}{a^{ny}-1}=\f
{a^x}{a^{ny}-1}\sum_{r=0}^{n-1}a^{ry}=\f{a^x}{a^y-1}\q\t{for}\q
n\in\Bbb Z^+.\tag 2.1$$
\par\q
\par{\bf Example 2.6.} Let $r>0$, $r\not=1$ and $\theta\in\Bbb R$, and let
$$f(x,y)=\f{r^{x+y}\cos(x-y)\theta-r^x\cos x\theta}{1-2r^y\cos
y\theta+r^{2y}}\q \t{and}\q g(x,y)=\f{r^{x+y}\sin(x-y)\theta-r^x\sin
x\theta}{1-2r^y\cos y\theta+r^{2y}}.$$ Then $f(x,y),g(x,y)\in I$.
\par By Euler's formula,
$$\align \f{(r\t{e}^{i\theta})^x}{(r\t{e}^{i\theta})^y-1}
&=\f{r^x(\cos x\theta+i\;\sin x\theta)}{r^y\cos y\theta-1+i\;\sin
y\theta\cdot r^y}=\f{r^x(\cos x\theta+i\;\sin x\theta)(r^y\cos
y\theta-1-i\;\sin y\theta\cdot r^y)}{(r^y\cos
y\theta-1)^2+r^{2y}\sin^2 y\theta}\\&=f(x,y)+ig(x,y).\endalign$$
Now, from Proposition 2.3 and (2.1) we deduce that $f(x,y),g(x,y)\in
I$.
\par{\bf Example 2.7.} Let $r>0$, $r\not=1$ and
$f(x,y)=\log \;(1-2r^{\f 1y}\cos\;2\pi\f xy+r^{\f 2y})$. Then $f\in
I$.
\par  For $n\in\Bbb Z^+$ and $z\in\Bbb C$,
$$ \prod_{j=0}^{n-1}(1-z\t{e}^{2\pi i\f jn})
=\prod_{j=0}^{n-1}(-\t{e}^{2\pi i\f jn})(z-\t{e}^{-2\pi i\f jn})
=(-1)^n\t{e}^{2\pi i\f{1+2+\cdots+n-1}n}(z^n-1)=1-z^n.$$ Set
$z=r^{\f 1{ny}}\t{e}^{2\pi i\f x{ny}}$. We get
$$\prod_{j=0}^{n-1}\big(1-r^{\f 1{ny}}\t{e}^{2\pi i\f {x+jy}{ny}}\big)
=1-r^{\f 1y}\t{e}^{2\pi i\f xy}.$$ Thus,
$$\sum_{j=0}^{n-1}\log \;|1-r^{\f 1{ny}}\t{e}^{2\pi i\f {x+jy}{ny}}|
=\log |1-r^{\f 1y}\t{e}^{2\pi i\f xy}|.$$ That is,
$$\sum_{j=0}^{n-1}\log \;\Big(1-2r^{\f 1{ny}}\cos\;
2\pi\f{x+jy}{ny}+r^{\f 2{ny}}\Big)=\log \Big(1-2r^{\f 1y}\cos\;
2\pi\f xy+r^{\f 2y}\Big).$$
\par{\bf Example 2.8.} Let $r>0$, $r\not=1$ and
 $f(x,y)=\f{r^{\f 1y}\sin 2\pi\f xy}{y(1-2r^{\f 1y}\cos\;2\pi\f
xy+r^{\f 2y})}$. Then $f\in I$.
\par This is immediate from Example 2.7 and Proposition 2.2.
\par{\bf Example 2.9.} Let $0<r<1$ and
 $f(x,y)=\f{1-r^{\f 2y}}{y(1-2r^{\f 1y}\cos\;2\pi\f
xy+r^{\f 2y})}$. Then $f\in I$.
\par Set
$$F(x,y)=\f 1{y(1-r^{\f 1y}\t{e}^{2\pi i\f xy})}-\f 1{2y}
=\f 1{2y}+\f 1y\sum_{k=1}^{\infty}r^{\f ky}\t{e}^{2\pi i\f{kx}y}.$$
By Example 2.1 and the proof of Proposition 2.7,
$\sum_{m=0}^{n-1}F(x+my,ny)=F(x,y)$ for $n\in\Bbb Z^+$. Note that
$$\align F(x,y)&=\f 1{y(1-r^{\f 1y}\cos 2\pi\f xy-i\;\sin 2\pi\f
xy\cdot r^{\f 1y})}-\f 1{2y}=\f{1-r^{\f 1y}\cos 2\pi\f xy+i\;\sin
2\pi\f xy\cdot r^{\f 1y}}{y(1-2r^{\f 1y}\cos\;2\pi\f xy+r^{\f
2y})}-\f 1{2y}
\\&=\f{1-r^{\f 2y}}{2y(1-2r^{\f 1y}\cos\;2\pi\f
xy+r^{\f 2y})}+i\f{r^{\f 1y}\sin 2\pi\f xy}{y(1-2r^{\f
1y}\cos\;2\pi\f xy+r^{\f 2y})}.\endalign$$ We see that $f(x,y)\in I$
by the above and Proposition 2.3.
\par{\bf Example 2.10.} For $y>0$ let
$$f(x,y)=\cases \log \big|2\;\sin\;\pi\f xy\big|&\t{if $\f xy\not\in\Bbb Z$,}
\\-\log y&\t{if $\f xy\in\Bbb Z$.}\endcases$$
Then $f\in I$.
\par Example 2.10 was essentially given in [7]. Here we give a
straightforward proof. For $x\not=2k\pi$ $(k\in\Bbb Z)$ we know that
$$\log \big|2\;\sin\f x2\big|=-\sum_{k=1}^{\infty}\f{\cos kx}k.$$
Thus, for $\f xy\notin\Bbb Z$ we have
$$\log \big|2\;\sin\pi\f xy\big|=-\f 1y\sum_{k=1}^{\infty}\f{\cos
2k\pi \f xy}{k/y}.$$ Since $\f xy\notin \Bbb Z$ implies
$\f{x+ry}{ny}\notin \Bbb Z$ for $n\in\Bbb Z^+$ and
$r\in\{0,1,\ldots,n-1\}$, from the above and Proposition 2.7 we
deduce that $\sum_{r=0}^{n-1}f(x+ry,ny)=f(x,y)$ for any $n\in\Bbb
Z^+$. Now assume that $\f xy\in\Bbb Z$. For $n\in\Bbb Z^+$ there is
a unique $r\in\{0,1,\ldots,n-1\}$ such that
$\f{x+ry}{ny}=\f{x/y+r}n\in\Bbb Z$. Hence,
$$\sum_{r=0}^{n-1}f(x+ry,ny)
=-\log ny+\sum\Sb r=0\\ n\nmid (x/y+r)\endSb^{n-1}\log\big|2\;\sin
\pi\f{x/y+r}n\big| =-\log ny+\sum_{k=1}^{n-1}\log 2\;\sin
\f{k\pi}n.$$ Since
$$\prod_{k=1}^{n-1}2\;\sin
\f{k\pi}n=\prod_{k=1}^{n-1}\f{\t{e}^{ik\pi/n}-\t{e}^{-ik\pi/n}}{i}
=\prod_{k=1}^{n-1}\big(1-\t{e}^{-\f{2k\pi i}n}\big)
=\lim_{x\rightarrow 1}\f{x^n-1}{x-1}=n,$$  we get
$$\sum_{r=0}^{n-1}f(x+ry,ny)=-\log ny+\log\prod_{k=1}^{n-1}2\;\sin
\f{k\pi}n =-\log ny+\log n=-\log y=f(x,y).$$

\par{\bf Example 2.11 ([7]).} For $y>0$ let
$$f(x,y)=\cases \f 1y\cot \pi\f xy&\t{if $\f xy\not\in\Bbb Z$,}
\\0&\t{if $\f xy\in\Bbb Z$.}\endcases$$
Then $f\in I$.

\par The Gamma function $\Gamma(x)$ is defined by
$$\Gamma(x)=\int_0^{\infty}t^{x-1}\hbox{e}^{-t}\;dt\q(x>0)\qtq{and}\Gamma(x+1)
=x\Gamma(x)\q(x\not=0,-1,-2,\ldots).$$
 For the properties of
$\Gamma(x)$ see [1] and [4].
\par{\bf Example 2.12 ([7]).} For $x\in\Bbb R$ and $y>0$ let
$$f(x,y)=\cases \log \Big|\f{y^{\f xy}\Gamma(\f xy)}{\sqrt{2\pi y}}
\Big|&\t{if $\f xy\not\in\{0,-1,-2,\ldots\}$,}
\\\log \f{y^{\f xy}}{(-\f xy)!}
\sqrt{2\pi y}&\t{if $\f xy\in\{0,-1,-2,\ldots\}$.}\endcases$$ Then
$f\in I$.
\par For $s>1$ and $x\not\in\{0,-1,-2,\ldots\}$ the Hurwitz zeta function
$\zeta(s,x)$ is defined by
$$\zeta(s,x)=\sum_{n=0}^{\infty}\f 1{(n+x)^s}.\tag 2.2$$
From [1, pp.51-52] we know that $\zeta(s,x)$ has a continuation to
the whole complex plane with a simple pole at $s=1$. That is,
$$\zeta(s,x)=\f{\Gamma(1-s)\t{e}^{-i\pi s}}{2\pi
i}\int_C\f{z^{s-1}\t{e}^{-xz}}{1-\t{e}^{-z}}dz\qtq{for}s\not=1,$$
where $C$ starts at infinity on the positive real axis, encircles
the origin once in the positive direction, excluding the points
$2k\pi i$ $(k\in\Bbb Z)$, and returns to positive infinity. In
particular, for $s<0$, Hurwitz showed that
$$\zeta(s,x)=\f{2\Gamma(1-s)}{(2\pi)^{1-s}}\Big(\sin\f{s\pi}2\sum_{k=1}^{\infty}
k^{s-1}\cos 2\pi kx+\cos\f{s\pi}2\sum_{k=1}^{\infty}k^{s-1}\sin 2\pi
kx\Big).\tag 2.3$$
 It is well-known that
$$\zeta(1-m,x)=-\f{B_m(x)}m\qtq{for}m\in\Bbb Z^+.\tag 2.4$$

\par{\bf Example 2.13.} For $s<0$, we have $y^{1-s}\zeta(s,\f
xy)\in I$ by (2.3) and Proposition 2.7. For $s>1$ define
$f(x,y)=y^{1-s}\zeta(s,\f xy)$ for $y>0$ and $\f
xy\not\in\{0,-1,-2,\ldots\}$. By (2.2) and Proposition 2.8,
$$\sum_{r=0}^{n-1}f(x+ry,ny)=f(x,y)\qtq{for any} n\in\Bbb Z^+.$$
\par{\bf Example 2.14.} Let
$$f(x,y)=\cases 1 &\t{if $\big\{\f xy\big\}<\f 12$,}
\\0&\t{if $\big\{\f xy\big\}=\f 12$,}
\\-1&\t{if $\big\{\f xy\big\}>\f 12$.}
\endcases$$
Then $f\in I$.

\section*{3. Main results for integrable invariant functions}
\par\q Suppose that $a(x),b(x)$ and $f(x,y)$ are real functions of
$x$ and their derivatives exist. The famous Leibniz's formula states
that
$$\f d{dx}\int_{a(x)}^{b(x)}f(x,t)dt=f(x,b(x))\f d{dx}b(x)-f(x,a(x))\f
d{dx}a(x)+\int_{a(x)}^{b(x)}\f{\partial}{\partial x}f(x,t)dt.\tag
3.1$$
 \pro{Theorem 3.1} Suppose that $f(x,y)$ is an integrable invariant function.
 \par $(\t{\rm i})$ We have
$$\int_x^{x+y}f(t,y)dt=\lim_{a\rightarrow 0^{+}}a f(x,a).$$
\par $(\t{\rm ii})$ We have
$$f(x+y,y)-f(x,y)= \f{d}
{dx}\lim_{a\rightarrow 0^{+}}af(x,a). $$ If $\f{\partial f}
{\partial x}$ is a piecewise continuous function of $x$, we also
have
$$f(x+y,y)-f(x,y)= \lim_{a\rightarrow 0^{+}}a\f{\partial f(x,a)}
{\partial x}.$$
\par $(\t{\rm iii})$ If $\f{\partial f} {\partial y}$ exists, then
$$f(x,y)=\int_x^{x-y}\f{\partial}{\partial y}f(t,y)\; dt.$$
\par $(\t{\rm iv})$ If $\f{\partial f} {\partial x}$ and
$\f{\partial f} {\partial y}$ exist, setting $g(x,y)=\f{\partial
f(x,y)} {\partial y}$ we have
$$\f{\partial f(x,y)} {\partial x}=g(x-y,y)-g(x,y).$$
\endpro
 {\it Proof}. Since
$\sum_{r=0}^{n-1}f(x+\f rn y,y)=f(x,\f yn)$ for $n\in\Bbb Z^+$, we
see that
$$\int_0^1f(x+uy,y)du=\lim_{n\rightarrow +\infty}
\f 1n\sum_{r=0}^{n-1}f\big(x+\f rny,y\big)
=\lim_{n\rightarrow +\infty}\f 1nf\big(x,\f yn\big)$$ and so
$$\int_x^{x+y}f(t,y)dt
=y\int_0^1f(x+uy,y)du=\lim_{n\rightarrow +\infty}\f ynf\big(x,\f
yn\big) =\lim_{a\rightarrow 0^{+}}a f(x,a).$$ This proves part(i).
By part(i) and (3.1), we obtain
$$f(x+y,y)-f(x,y)=\f{\partial
} {\partial x}\int_x^{x+y}f(t,y)dt= \f{d} {dx}\lim_{a\rightarrow
0^{+}}af(x,a).$$ If $\f{\partial f} {\partial x}$ is a piecewise
continuous function of $x$ and $h(x,y)=\f{\partial f(x,y)} {\partial
x}$, we see that
$$f(x+y,y)-f(x,y)=\int_x^{x+y}h(t,y)dt=\lim_{a\rightarrow
0^{+}}ah(x,a).$$ This proves part (ii).

\par Suppose that $\f{\partial f}
{\partial y}$ exists and $g(x,y)=\f{\partial}{\partial y}f(x,y)$. By
part (i),
$$\f{\partial}{\partial y}\int_x^{x+y}f(t,y)dt
= \f{d}{dy} \lim_{a\rightarrow 0^{+}} a f(x,a)=0.$$ By (3.1),

$$\f{\partial}{\partial
y}\int_x^{x+y}f(t,y)dt=f(x+y,y)+\int_x^{x+y}g(t,y)dt.$$  Thus,
$$f(x+y,y)=-\int_x^{x+y}g(t,y)dt\qtq{and so}f(x,y)=\int_x^{x-y}g(t,y)\; dt.$$
This proves part (iii).

\par Now suppose that $\f{\partial f} {\partial x}$ and
$\f{\partial f} {\partial y}$ exist and $g(x,y)=\f{\partial f(x,y)}
{\partial y}$. By (3.1) and part (iii), we obtain part (iv).
\par Summarizing the above proves the theorem.
\par\q

\pro{Theorem 3.2} Suppose that $f\in I$ and $g(x,y)=\f{\partial
f(x,y)}{\partial y}$.
\par $(\t{\rm i})$ If $g(x,y)$ is an even function of
$x$, then
$$f(y-x,y)=f(x,y)\qtq{and}\int_0^{\f y2}f(t,y)dt
=\f 12\lim_{a\rightarrow 0^{+}} a f(0,a).$$
\par $(\t{\rm ii})$ If $g(x,y)$ is an odd function of
$x$, then
$$f(y-x,y)=-f(x,y)\qtq{and}\int_0^yf(t,y)dt
=\lim_{a\rightarrow 0^{+}} a f(0,a)=0.$$
\endpro
{\it Proof.} If $g(-x,y)=(-1)^m g(x,y)$, from Theorem 3.1(iii) we
see that
$$ f(y-x,y)=\int_{y-x}^{-x}g(t,y)dt=\int_x^{x-y}g(-u,y)du
=(-1)^m\int_x^{x-y}g(u,y)du=(-1)^mf(x,y).$$ Hence,
$$\int_{\f y2}^yf(t,y)dt=\int_0^{\f y2}f(y-t,y)dt=(-1)^m\int_0^{\f
y2}f(t,y)dt$$ and so $\int_0^yf(t,y)dt=(1+(-1)^m)\int_0^{\f
y2}f(t,y)dt.$ This together with Theorem 3.1(i) (with $x=0$) yields
the result.

\par\q
\par{\bf Remark 3.1} From Theorems 3.1-3.2 and Examples 2.7, 2.10 and
2.12 (with $y=1$) we may deduce the following known integrals:
$$\align &(\t{Euler})\q \int_0^{\f {\pi}2}\log \sin\; xdx=-\f
{\pi}2\log 2,\tag 3.2
\\&(\t{Poission})\q \int_0^{\pi}\log(1-2r\cos x+r^2)dx=\cases
2\pi\log r&\t{if $r>1$},\\0&\t{if $0<r<1$,}\endcases\tag 3.3
\\&(\t{Raabe})\q \int_a^{a+1}\log\Gamma(x)dx=a(\log
a-1)+\log\sqrt{2\pi}\ (a>0). \tag 3.4\endalign$$

 \pro{Theorem 3.3} Suppose that $f\in I$ and $\f{\partial}{\partial
y}f(x,y)$ exists. For fixed $y>0$, $f(x,y)$ is a function of $x$
with bounded variation at the closed interval $[a,b]$.\endpro
 {\it Proof}.
Suppose $g(x,y)=\f{\partial}{\partial y}f(x,y)$. Then
$f(x,y)=\int_x^{x-y}g(t,y)dt$ by Theorem 3.1(iii). For
$a=x_0<x_1<\cdots <x_n=b$ we see that
$$\align\sum_{i=0}^{n-1}|f(x_{i+1})-f(x_i)|
&=\sum_{i=0}^{n-1}\Big|\int_{x_{i+1}}^{x_{i+1}-y}g(t,y)dt-
\int_{x_{i}}^{x_{i}-y}g(t,y)dt\Big|
\\&=\sum_{i=0}^{n-1}\Big|\int_{x_{i+1}}^{x_i}g(t,y)dt+
\int_{x_i-y}^{x_{i+1}-y}g(t,y)dt\Big|
\\&\le \sum_{i=0}^{n-1}\Big|\int_{x_i}^{x_{i+1}}g(t,y)dt\Big|
+\sum_{i=0}^{n-1}\Big|\int_{x_i-y}^{x_{i+1}-y}g(t,y)dt\Big|
\\&\le  \sum_{i=0}^{n-1}\int_{x_i}^{x_{i+1}}|g(t,y)|dt
+\sum_{i=0}^{n-1}\int_{x_i-y}^{x_{i+1}-y}|g(t,y)|dt
\\&=\int_a^b|g(t,y)|dt+\int_{a-y}^{b-y}|g(t,y)|dt
.\endalign$$ This proves the theorem.
\par\q
 \pro{Theorem 3.4 (Product Theorem) } For $g,h\in I$ define
$$g*h(x,y)=\int_0^xg(t,y)h(x-t,y)dt+\int_x^yg(t,y)h(x+y-t,y)dt.$$
Then $g*h\in I$ and $$\int_0^yg*h(x,y)dx=\int_0^yg(x,y)dx\cdot
\int_0^yh(x,y)dx.$$
\endpro
\par {\it Proof.} Set $f(x,y)=g*h(x,y)$,
$h_1(x,y)=h(y-x,y)$ and
$$f_1(x,y)=\int_0^yh_1\Big(y\big\{\f{t+x}y\big\},y\Big)g(t,y)dt.$$
Since $h(x,y)\in I$, we have $h_1(x,y)\in I$ and
$h_1(y\{\f{t+x}y\},y)\in I$ by Propositions 2.5 and 2.6. Hence, for
$n\in\Bbb Z^+$ we have
$$\align \sum_{r=0}^{n-1}f_1(x+ry,ny)&=\sum_{r=0}^{n-1}\int_0^{ny}
h_1\Big(ny\big\{\f{t+x+ry}{ny}\big\},ny\Big)g(t,ny)dt
\\&=\int_0^{ny}h_1\Big(y\big\{\f{t+x}y\big\},y\Big)g(t,ny)dt\\
&=\sum_{r=0}^{n-1}\int_{ry}^{(r+1)y}h_1\Big(y\big\{\f{t+x}y\big\},y\Big)
g(t,ny)dt\\
&=\sum_{r=0}^{n-1}\int_0^yh_1\Big(y\big\{\f{u+ry+x}y\big\},y\Big)g(u+ry,ny)du\\
&=\int_0^yh_1\Big(y\big\{\f{u+x}y\big\},y\Big)\sum_{r=0}^{n-1}g(u+ry,ny)du\\
&=\int_0^yh_1\Big(y\big\{\f{u+x}y\big\},y\Big)g(u,y)du=f_1(x,y).
\endalign$$
Thus $f_1\in I$. For $0\le x\le y$ we see that
$$\align f_1(x,y)&=\int_0^yh_1\big(y\big\{\f{t+x}y\big\},y\big)g(t,y)dt
\\&=\int_0^{y-x}h_1\big(y\cdot \f{t+x}y,y\big)g(t,y)dt+\int_{y-x}^yh_1\big(y
\big(\f{t+x}y-1\big),y\big)g(t,y)dt\\
&=\int_0^{y-x}h_1(t+x,y)g(t,y)dt+\int_{y-x}^yh_1(t+x-y,y)g(t,y)dt
\\&=\int_0^{y-x}h(y-x-t,y)g(t,y)dt+\int_{y-x}^yh(2y-x-t)g(t,y)dt
\\&=g*h(y-x,y)=f(y-x,y)\endalign$$
and so $f(x,y)=f_1(y-x,y)\in I$ by Proposition 2.5.
\par Set ${\bar F}(x,y)=F(x+y,y)-F(x,y)$. Then
$$\align {\bar
f}(x,y)&=f(x+y,y)-f(x,y)\\&=\int_0^{x+y}g(t,y)h(x+y-t,y)dt+\int_{x+y}^yg(t,y)h(x+2y-t,y)dt
\\&\q-\int_0^xg(t,y)h(x-t,y)dt-\int_x^yg(t,y)h(x+y-t,y)dt
\\&=\int_0^xg(t,y)(h(x-t+y,y)-h(x-t,y))dt\\&\q+\int_y^{x+y}(h(x+y-t,y)
-h(x+2y-t,y))g(t,y)dt
\\&=\int_0^x g(t,y){\bar
h}(x-t,y)dt+\int_0^x(h(x-t,y)-h(x+y-t,y))g(t+y,y)dt
\\&=-\int_0^x{\bar h}(x-t,y){\bar g}(t,y)dt.
\endalign$$
For $F(x,y)\in I$ and $n\in\Bbb Z^+$, we see that
$$\align {\bar F}(x,ny)&=F(x+ny,ny)-F(x,ny)
=\sum_{r=0}^{n-1}F(x+y+ry,ny)-\sum_{r=0}^{n-1}F(x+ry,ny)
\\&=F(x+y,y)-F(x,y)={\bar F}(x,y).\endalign$$
Hence, for any $x\in\Bbb R$ and $n\in\Bbb Z^+$,
$${\bar f}(x,ny)=-\int_0^x{\bar h}(x-t,ny){\bar g}(t,ny)dt
=-\int_0^x{\bar h}(x-t,y){\bar g}(t,y)dt={\bar f}(x,y).\tag 3.5$$
Therefore, if $\sum_{r=0}^{n-1}f(x+ry,ny)=f(x,y)$, then
$$\align \sum_{r=0}^{n-1}f(x+y+ry,ny)
&=\sum_{r=0}^{n-1}f(x+ry,ny)+f(x+ny,ny)-f(x,ny)
\\&=f(x,y)+f(x+y,y)-f(x,y)=f(x+y,y)\endalign$$
and
$$\align \sum_{r=0}^{n-1}f(x-y+ry,ny)
&=\sum_{r=0}^{n-1}f(x+ry,ny)-f(x-y+ny,ny)+f(x-y,ny)
\\&=f(x,y)-(f(x-y+ny,ny)-f(x-y,ny))
\\&=f(x,y)-(f(x-y+y,y)-f(x-y,y)) =f(x-y,y).\endalign$$
Since we have proved that $\sum_{r=0}^{n-1}f(x+ry,ny)=f(x,y)$ for
$0\le x\le y$, we must have $\sum_{r=0}^{n-1}f(x+ry,ny)=f(x,y)$ for
any $x\in\Bbb R$ and $y>0$. Hence $f(x,y)\in I$.

\par For fixed $y>0$ let $H(x,y)$ be a primitive function of $h(x,y)$.
That is, $\f{\partial\; H}{\partial x}=h$. Set
$$g\ast H(x,y)=\int_0^xH(x-t,y)g(t,y)dt+\int_x^yH(x+y-t,y)g(t,y)dt.$$
Then
$$g\ast H(y,y)=\int_0^yH(y-t,y)g(t,y)dt=g\ast H(0,y).$$
Appealing to (3.1) we see that
$$\align \f {\partial}{\partial x}g\ast
H(x,y)&=\int_0^xh(x-t,y)g(t,y)dt+\int_x^yh(x+y-t,y)g(t,y)dt
\\&+H(0,y)g(x,y)
-H(y,y)g(x,y)\\
&=g*h(x,y)-g(x,y)\int_0^yh(x,y)dx
\endalign$$
and so
$$\align 0&=g\ast H(y,y)-g\ast H(0,y)=\int_0^y\f {\partial}{\partial
x}g\ast
H(x,y)dx\\
&=\int_0^yg*h(x,y)dx-\int_0^yh(x,y)dx\int_0^yg(x,y)dx.
\endalign$$
This completes the proof.
\par\q

\pro{Corollary 3.1} Suppose that $g(x,y)\in I$, $a>0$, $a\not=1$ and
$$f(x,y)=\f{a^x}{a^y-1}\int_0^ya^{-t}g(t,y)dt+a^x\int_x^ya^{-t}g(t,y)dt.$$
Then $f(x,y)\in I$.
\endpro
{\it Proof.} Since $\f{a^x}{a^y-1}\in I$ by Example 2.5 and
$$\align\f{a^x}{a^y-1}*g(x,y)
&=\int_0^x\f{a^{x-t}}{a^y-1}g(t,y)dt+\int_x^y\f{a^{x+y-t}}{a^y-1}g(t,y)dt
\\&=\int_0^x\f{a^{x-t}}{a^y-1}g(t,y)dt+\int_x^y\f{a^{x-t}}{a^y-1}g(t,y)dt
+\int_x^ya^{x-t}g(t,y)dt =f(x,y),\endalign$$  the result follows
from Theorem 3.4.

 \pro{Theorem 3.5} Suppose $f\in I$ and
$$F(x,y)=\int_y^xf(t,y)dt+\f 1y\int_0^ytf(t,y)dt.$$
Then $F\in I$ and $\f{\partial}{\partial x}F(x,y)=f(x,y)$.
\endpro
{\it Proof.} By Example 2.2 (with $m=1$), $\f xy-\f 12\in I$. Set
$$G(x,y)=\int_0^x\Big(\f{x-t}y-\f 12\Big)f(t,y)dt+\int_x^y
\Big(\f{x+y-t}y-\f 12\Big)f(t,y)dt.$$ Then $G\in I$ by Theorem 3.4.
It is clear that
$$\align G(x,y)=\int_0^y\Big(\f{x-t}y-\f 12\Big)f(t,y)dt+\int_x^yf(t,y)dt
=\Big(\f xy-\f 12\Big)\int_0^yf(t,y)dt-F(x,y).
\endalign$$
By Theorem 3.1, $\int_0^yf(t,y)dt$ is a constant and so $(\f xy-\f
12)\int_0^yf(t,y)dt\in I$. Hence $F(x,y)\in I$. By Leibniz's formula
(3.1), $\f{\partial}{\partial x}F(x,y)=f(x,y)$. Thus the theorem is
proved.
\par\q
\par
For $m,n\in\Bbb Z^+$, from [4, pp.26-27] we have
$$\align &B_{2m+1}=0, \q
B_m(1-t)=(-1)^mB_m(t),\q\f{d}{dx}B_m(x)=mB_{m-1}(x),\tag 3.6
\\&\int_0^1B_n(t)dt=\f{B_{n+1}(1)-B_{n+1}}{n+1}
=((-1)^{n+1}-1)\f{B_{n+1}}{n+1}=0,\tag 3.7
\\&\int_0^1B_m(t)B_n(t)dt=(-1)^{m-1}\f{B_{m+n}}{\b{m+n}m}.\tag 3.8
\endalign$$
Thus,
$$\int_0^1B_m(1-t)B_n(t)dt=(-1)^m\int_0^1B_m(t)B_n(t)dt =-
\f{B_{m+n}}{\b{m+n}m}=- \f{B_{m+n}(1)}{\b{m+n}m}.\tag 3.9$$

\pro{Theorem 3.6} For any positive integers $m$ and $n$ we have
$$B_{m+n}(x)=-\binom{m+n}m\Big(\int_0^1B_m(x-t)B_n(t)dt
+m\int_x^1(x-t)^{m-1}B_n(t)dt\Big)\tag 3.10$$ and
$$\f{-y^{m-1}B_m\sls xy}{m!}*\f{-y^{n-1}B_n\sls xy}{n!}=
\f{-y^{m+n-1}B_{m+n}\sls xy}{(m+n)!}.$$
\endpro
{\it Proof.}  We first prove (3.10) by induction on $m$. Since
$B_1(x)=x-\f 12$, using (3.1), (3.6) and (3.7) we see that
$$\align&\f d{dx}\Big(\int_0^1B_1(x-t)B_n(t)dt+\int_x^1B_n(t)dt\Big)
\\&=\int_0^1B_n(t)dt-B_n(x)=-B_n(x)
=\f d{dx}\Big(-\f {B_{n+1}(x)}{n+1}\Big).\endalign$$ By (3.9),
$\int_0^1B_1(1-t)B_n(t)dt=- \f{B_{n+1}(1)}{n+1}$. Therefore,
$$\int_0^1B_1(x-t)B_n(t)dt+\int_x^1B_n(t)dt=-\f {B_{n+1}(x)}{n+1}.$$
This shows that the result is true for $m=1$.
\par Now assume that (3.10) holds for $m=k$ $(k\in\Bbb Z^+)$.
Appealing to (3.1) and (3.6), we see that
$$\align &\f d{dx}\Big(\int_0^1B_{k+1}(x-t)B_n(t)dt
+(k+1)\int_x^1(x-t)^kB_n(t)dt\Big)
\\&=\int_0^1(k+1)B_k(x-t)B_n(t)dt+(k+1)k\int_x^1(x-t)^{k-1}B_n(t)dt
\\&=-\f {k+1}{\b{k+n}{k}}B_{k+n}(x)
=\f d{dx}\Big(-\f {B_{k+1+n}(x)}{\b{k+1+n}{k+1}}\Big).
\endalign$$
 By (3.9),
$$\int_0^1B_{k+1}(1-t)B_n(t)dt
+(k+1)\int_1^1(1-t)^kB_n(t)dt =- \f{B_{k+1+n}(1)}{\b{k+1+n}{k+1}}.$$
Therefore,
$$-\f {B_{k+1+n}(x)}{\b{k+1+n}{k+1}}=\int_0^1B_{k+1}(x-t)B_n(t)dt
+(k+1)\int_x^1(x-t)^kB_n(t)dt.$$ This shows that (3.10) is true for
$m=k+1$. Hence (3.10) is proved by induction.
\par From [4] we know that $B_m(x+1)=B_m(x)+mx^{m-1}$.
Thus, using (3.10) we deduce that
$$\align &-y^{m-1}\f{B_m(\f xy)}{m!}*\Big(-y^{n-1}\f{B_n\sls
xy}{n!}\Big)
\\&=\int_0^xy^{n-1}\f{B_n\sls ty}{n!}\cdot
y^{m-1}\f{B_m\sls{x-t}y}{m!} dt+\int_x^yy^{n-1}\f{B_n\sls ty}{n!}
\cdot y^{m-1}\f{B_m\sls{x+y-t}y}{m!}dt
\\&=\f{y^{m+n-2}}{m!n!}\Big(\int_0^yB_n\Ls tyB_m\Ls{x-t}ydt
+\int_x^yB_n\Ls ty\Big(B_m\Ls{x+y-t}y-B_m\Ls{x-t}y\Big)dt\Big)
\\&=\f{y^{m+n-2}}{m!n!}\Big(\int_0^yB_m\Ls{x-t}yB_n\Ls tydt
+\int_x^y m\Ls {x-t}y^{m-1}B_n\Ls tydt\Big)
\\&=\f{y^{m+n-1}}{m!n!}\Big(\int_0^1B_m\Big(\f xy-u\Big)B_n(u)du+\int_{\f
xy}^1 m\Big(\f xy-u\Big)^{m-1}B_n(u)du\Big)
\\&=\f{y^{m+n-1}}{m!n!}\Big(-\f{B_{m+n}\sls xy}{\b{m+n}m}\Big)
=-y^{m+n-1}\f{B_{m+n}\sls xy}{(m+n)!}.
\endalign$$
This completes the proof.
\par{\bf Remark 3.2} For $\alpha,\beta>1$, from (2.3) we may deduce
that
$$\f{\zeta(1-\alpha,\f xy)y^{\alpha-1}}{\Gamma(\alpha)}
*\f{\zeta(1-\beta,\f xy)y^{\beta-1}}{\Gamma(\beta)}
=\f{\zeta(1-\alpha-\beta,\f
xy)y^{\alpha+\beta-1}}{\Gamma(\alpha+\beta)}.
$$
This is the generalization of Theorem 3.6. We omit its proof, which
is somewhat complicated.

\section*{Acknowledgements} The author
thanks Zhi-Wei Sun for his helpful comments on the proof of Theorem
3.4, and the author was supported by the National Natural Science
Foundation of China (grant No. 12271200).


\begin{thebibliography}{9}
\bibitem [1] {} G.E. Andrews, R. Askey and R. Roy, Special
Functions, Cambridge Univ. Press, Cambridge, 1999.
\bibitem [2] {} H. Bass, Generators and relations for cyclotomic units, Nagoya
Math. J. 27(1966), 401-407.
\bibitem [3] {} D.S. Kubert, The universal ordinary distribution, Bull. Soc.
Math. France 107(1979), 179-202.
\bibitem [4] {}
W. Magnus, F. Oberhettinger and R.P. Soni, Formulas and Theorems for
the Special Functions of Mathematical Physics (3rd edit.), Springer,
New York, 1966, pp. 25-27.
\bibitem [5] {} J. Milnor, On polylogarithms, Hurwitz zeta functions, and the
Kubert identities, Enseign. Math. 29(1983), 281-322.
\bibitem [6] {}
Z.W. Sun, Systems of congruences with multipliers, Nanjing Univ. J.
-Math. Biq. 6(1989), 124-133.
\bibitem [7] {} Z.W. Sun, On covering invariant, in: Analytic Number
Theory (Beijing/Kyoto,1999), 277-302, Dev. Math., 6, Kluwer Acad.
Publ., Dordrecht, 2002.
\bibitem [8] {} Z.W. Sun, Arithmetic properties of periodic maps,
Math. Res. Lett. 11(2004), 187-196.
\bibitem [9] {} H. Walum, Multiplication formulae for periodic
functions, Pacific J. Math. 149(1991), 383-396.


\end{thebibliography}
\end{document}